\documentclass[12pt]{article}
\usepackage{latexsym,amsfonts,amssymb}
\usepackage[dvips]{color}
\usepackage{colordvi,multicol}

\def \cal{\mathcal}

\newtheorem{thm}{Theorem}[section]
\newtheorem{cor}[thm]{Corollary}
\newtheorem{lem}[thm]{Lemma}

\newtheorem{defi}[thm]{Definition}
\newtheorem{rem}[thm]{Remark}

\begin{document}

\title{\bf On variants of conflict-free-coloring for hypergraphs}
\author{Zhen Cui and Ze-Chun Hu\thanks{Corresponding
author: Department of Mathematics, Nanjing University, Nanjing
210093, PR China\vskip 0cm E-mail address: huzc@nju.edu.cn}\\
 {\small Department of Mathematics, Nanjing University}}
 \date{}
 \maketitle

  \vskip 1.4cm

\noindent{\bf Abstract}\quad Conflict-free coloring is a kind of
vertex coloring of hypergraphs requiring each hyperedge to have a
color which appears only on one vertex. More generally,  for a positive
integer $k$ there are $k$-conflict-free colorings ($k$-CF-colorings
for short) and $k$-strong-conflict-free colorings ($k$-SCF-colorings
for short). 
Let $H_n$ be the
hypergraph of which the vertex-set is $V_n=\{1,2,\dots,n\}$ and the
hyperedge-set $\cal{E}_n$ is the set of all (non-empty) subsets of
$V_n$ consisting of consecutive elements of $V_n$. Firstly, we study
the $k$-SCF-coloring of $H_n$, give the exact $k$-SCF-coloring
chromatic number of $H_n$ for $k=2,3$, and present upper and
lower bounds of the $k$-SCF-coloring chromatic number of $H_n$ for all $k$.
Secondly, we give the exact $k$-CF-coloring chromatic number of $H_n$ for all
$k$.


%


\section{Introduction}

A {\it hypergraph} is a pair $H=(V,\mathcal E)$ where $V$ is a set and
$\mathcal E$ is a collection of subsets of ${V}$. The elements of
$V$ are called {\it vertices} and the elements of $\mathcal E$ are
called {\it hyperedges}. If for any $e \in \mathcal E$, $|e|=2$,
then the pair $(V,\mathcal E)$ is a $simple \ graph$. For a subset
$V^{'} \subset V$, we call the hypergraph $H(V^{'})=(V^{'},\{S \cap
V^{'}\vert S \in \mathcal E \})$  the {\it sub-hypergraph} induced
by $V^{'}$.  An $m$-{\it  coloring} for some $m \in \mathbb N$ of
(the vertices of) $H$ is a function $ \phi : V \to \{ 1,\ldots ,m
\}$. Let $ \phi $ be an $m$-{\it coloring} of $H$, if for any $e \in
\mathcal E$ with $|e|\geq 2$, there exist at least two vertices $x,y
\in e$ such that $ \phi(x) \neq \phi(y)$, we call $\phi$ {\it
proper} or {\it non-monochromatic}. Let $\chi(H)$ denote the least
integer $m$ for which $H$ admits a proper coloring with $m$ colors.
The following coloring is more restrictive than non-monochromatic
coloring.

\begin{defi}[Conflict-Free Coloring]
Let $H=(V,\mathcal E)$ be a hypergraph and let $\phi:V \to \{
1,\ldots,m \}$ be some coloring of $H$.  $\phi$ is called a
conflict-free coloring (CF-coloring for short) if for any $e \in
\mathcal E $ there is a vertex $x \in e$ such that $\forall y \in e,
y \neq x \Rightarrow \phi(y) \neq \phi(x)$.
\end{defi}
The CF-coloring chromatic number $\chi_{cf}(H)$ is the least integer $m$ for which $H$ admits
a CF-coloring with $m$ colors.

The notion of CF-coloring was first introduced and studied by
Smorodinsky \cite{38} and  Even et al. \cite{18}. Such coloring is
very useful in wireless networks, radio frequency identification
(RFID) networks and vertex ranking prolem. Refer to the survey paper
\cite{Sm11} by Smorodinsky and the references therein for more
descriptions. Such coloring have attracted many researchers from the
computer science and mathematics community. As to CF-coloring of
hypergraphs that arise in geometry, refer to Smorodinsky \cite{38},
Even et al. \cite{18}, Har-Peled and Smorodinsky \cite{20},
Smorodinsky \cite{39}, Pach and Tardos \cite{33}, Ajwani et al.
\cite{2}, Chen et al. \cite{13},
 Alon and Smorodinsky \cite{3},
Lev-Tov and
Peleg \cite{28} and etc. As to
CF-coloring of arbitrary hypergraphs, refer to
Pach and Tardos
\cite{32}.

Smorodinsky  \cite{38} considered extensions of CF-coloring and
introduced the following notion.

\begin{defi}
[$k$-CF-coloring]
 Let $H=(V,\mathcal E)$ be a hypergraph, $k$ be a positive integer. A coloring
$\phi:V\to\{1,\ldots,m\}$ is
 called a $k$-CF-coloring of $H$ if
 for any $e \in \mathcal E$ there is a color $j$ such that $1\leq
 |\{ v\in e |\phi(v) =j \}| \leq k$.
\end{defi}
The $k$-CF-coloring chromatic number $\chi_{kcf}(H)$ is the smallest number of
colors in any possible $k$-CF-coloring of $H$. Note that
1-CF-coloring of a hypergraph is simply a CF-coloring.

Refer to Smorodinsky \cite{38} and Har-Peled and Smorodinsky
\cite{20} for the study of $k$-CF-coloring. Another extension of
CF-coloring is called $k$-SCF-coloring, which is defined as follows:

\begin{defi}
[$k$-SCF-coloring] Let $H=(V,\mathcal{E})$ be a hypergraph, $k$ be a
positive integer. A coloring $\phi:V\to\{1,\ldots,m\}$ is called a
 $k$-SCF-coloring
if for any $e\in \mathcal E$ with $|e| \geq k$, there are at least
$k$ colors which appear only once in $e$, and for any $e\in \mathcal
E$ with $|e| < k$ all points in $e$ are of different colors.
\end{defi}
The $k$-SCF-coloring chromatic number $\chi_{kscf}(H)$ is the smallest number of colors in any possible
$k$-SCF-coloring of $H$. Note that 1-SCF-coloring is just a
CF-coloring.

Abellanas et al. \cite{1} were the first to study
$k$-SCF-coloring\footnote{They referred to such a coloring as
$k$-conflict-free coloring}. They focused on the special case of hypergraphs induced by $n$ points in $\mathbb{R}^2$ where hyperedges are cutdown by discs (which means that a hyperedge is composed with all points that are in a disc) and showed that in this case the hypergraph admits a $k$-SCF-coloring with $O(\frac{\log n}{\log\frac{ck}{ck-1}})=O(k\log n)$ colors for some absolute constant $c$. Aloupis et al. \cite{5} introduced
another coloring called $k$-colorful coloring, which has interesting
connection with strong-conflict-free coloring. Refer to Horev et al.
\cite{22} for the connection and research for $k$-SCF-coloring.

Throughout the rest of this paper, we let $H_n=(V_n,\mathcal{E}_n)$
be the complete hypergraph over $n$ points, here
$V_n=\{1,2,\dots,n\}$, and $\mathcal{E}_n$ is the set of all (non-empty) subsets of $V_n$ consisting of
consecutive elements of $V_n$.  For example, for $n=4$, we have $V_4=\{1,2,3,4\}$ and
$$
\mathcal{E}_4=\{\{1\},\{2\},\{3\},\{4\},\{1,2\},\{2,3\},\{3,4\},\{1,2,3\},\{2,3,4\},\{1,2,3,4\}\}.
$$
Har-Peled and Smorodinsky  \cite{20} proves that
$\chi_{cf}(H_n)=\lfloor\log n\rfloor+1$ as a simple yet an important
example of CF-coloring of a hypergraph.

In Section 2, we consider $k$-SCF-coloring of $H_n$, give the
exact $k$-SCF-coloring chromatic number of $H_n$ for $k=2,3$, and present upper and lower bounds for the $k$-SCF-coloring chromatic number of $H_n$ for all $k$.  In Section 3, we give the exact $k$-CF-coloring chromatic number of $H_n$ for all $k$.




\section{$k$-SCF-coloring of $H_n$}\label{kSCF}
In this section, we consider $k$-SCF-coloring of $H_n$ and focus on
$\chi_{kscf}(H_n)$ especially. By the definition of $k$-SCF-coloring, it's obvious that any $k$-SCF-coloring of hypergraphs induces a $k$-SCF-coloring on their sub-hypergraphs, and so $\chi_{kscf}(H_n)$ is monotonic in $n$. For $m\geq 1,$ define
\begin{equation}
g_k(m)=\sup\{n: \chi_{kscf}(H_n)\leq m\},
\end{equation}
i.e. $g_k(m)$ is the largest $n$ such that  we can give $H_n$ a
$k$-SCF-coloring by using $m$ colors. The idea is that  if we can
get one clear expression of $g_{k}(m)$ as a function of $m$, then we
will be able to obtain $\chi_{kscf}(H_n)$ by the formula
$\chi_{kscf}(H_n)=\inf\{m:g_k(m)\geq n\}$.
 Generally, we have
the following inequalities.

\begin{lem}\label{lem2.1}
Suppose $k,m,p,q\in\mathbb{N}$,  $m\geq k$ and $k+1=p+q$. Then we
have
$$g_{k}(m)\ {\leq}\ g_{k}(m-p)+g_{k}(m-q)+1.$$
\end{lem}
\noindent {\bf Proof.} Suppose the inequality is not true, then
there is some way to color $g_{k}(m-p)+g_{k}(m-q)+2$ points using
$m$ colors and the coloring is $k$-SCF. Suppose these
$g_{k}(m-p)+g_{k}(m-q)+2$ points are
$$\underbrace{\scriptstyle{1,\ 2,\ \dots,\ g_{k}(m-p)-1,\
g_{k}(m-p),}}_{\textstyle{M\ \textit{region}}}\ A,\
\underbrace{\scriptstyle{g_{k}(m-p)+2,\ g_{k}(m-p)+3,\ \dots,\
g_{k}(m-p)+g_{k}(m-q)+1,}}_{\textstyle{N\ \textit{region}}}\ B,$$
where $A,B$ are two points. Because the coloring is $k$-SCF, there
are $k$ colors $\{a_{1},\dots,a_{k}\}$, which appear only once over
the $g_{k}(m-p)+g_{k}(m-q)+2$ points. But less than $p$  of these
$k$ colors could appear in $\{N\ \textit{region}\}\cup\{B\}$,
otherwise
 there will be at most $m-p$ colors in $\{M\ \textit{region}\}\cup\{A\}$, which is a contradiction with respect to (w.r.t. for short) the definition of $g_{k}(m-p)$. So there are at most $p-1$ colors of $\{a_{1},\dots,a_{k}\}$ which appear in $\{N\ \textit{region}\}\cup\{B\}$, that is, at least $k-(p-1)=q$ colors of $\{a_{1},\dots,a_{k}\}$ appear in $\{M\ \textit{region}\}\cup\{A\}$. So now there will be at most $m-q$ colors which appear in $\{N\ \textit{region}\}\cup\{B\}$. This  is a contradiction w.r.t. the definition of $g_{k}(m-q)$. Hence we have that $g_{k}(m){\leq}g_{k}(m-p)+g_{k}(m-q)+1$.
\hfill\fbox

By the above lemma, we have the following result.
\begin{cor}\label{cor2.2}
Suppose $k,m,p\in\mathbb{N}$,  $m\geq k$. Then we have\\
(i) if $k=2p$, then $g_{k}(m)\leq g_{k}(m-p)+g_{k}(m-p-1)+1$;\\
(ii) if $k=2p+1$, then $g_{k}(m)\leq 2g_{k}(m-p-1)+1$.
\end{cor}

Next we focus on  two simple cases $k=2$ and $k=3$.

\begin{thm}\label{thm2.2}
For any integer $m\geq 2$, we have
\begin{eqnarray}\label{thm2.2-1}
g_2(m+1)=g_2(m)+g_2(m-1)+1.
\end{eqnarray}
\end{thm}
{\bf Proof.} By Corollary \ref{cor2.2}(i), in order to prove
(\ref{thm2.2-1}), we  need  only to prove  that for any $m\geq 2,$
$g_2(m+1)\geq g_2(m)+g_2(m-1)+1$, i.e. there exists  a
2-SCF-coloring for $g_{2}(m)+g_{2}(m-1)+1$ points using $m+1$
colors. For simplicity, in the folloing we use  $(a_1,a_2,\ldots,a_n)$ to denote a
coloring $\phi:\{1,2,\ldots,n\}\to \{1,2,\ldots,m\}$ for some
$m$ with  $a_i=\phi(i)$ for each $i=1,\ldots,n$, and denote by $C_m$ one
2-SCF-coloring for $g_{2}(m)$ points using $m$ colors. If $C_m=(a_1,a_2,\ldots,a_{n-1},a_n)$, then we denote
$C_m^{-1}=(a_n,a_{n-1},\ldots,a_2,a_1)$, and it's trivivally a 2-SCF-coloring.

Obviously, we
have $g_2(1)=1,g_2(2)=2$. Now we give the constructions by the
inductive method. For $m=1,2,3$ and 4, let
\begin{eqnarray*}
C_1=({\bf 1}),\ C_2=({\bf 1},{\bf 2}),\ C_3=(1,{\bf 2},{\bf 3},1),\ C_4=(1,2,{\bf 3},1,{\bf 4},2,1), \end{eqnarray*}
where the bold colors appear only once. We can easily check  that all of  them are 2-SCF-coloring.

For $m=5$, we construct a coloring  as follows:
\begin{eqnarray*}
(1,2,3,1,{\bf 4},2,1,{\bf 5},2,3,1,2),
\end{eqnarray*}
which is composed of several parts, denoted by
\begin{eqnarray*}
C_5=(\underbrace{C_{3},\ \  {\bf 4},\ \  M_{5}}_{C_4},\ \  {\bf 5},\ \  T_{5}).
\end{eqnarray*}
The middle part $M_{5}=(2,1)$ has a special property, namely when it is read from right to left, denoted $M_{5}^{-1}=(1,2)$, it could be derived from $M_5$ by substituting 1 with 2, and 2 with 1. We write $M_{5}^{-1}=P(1,2)M_{5}$, where $P(1,2)$ stands for the permutation exchanging 1 and 2.
Next, for the tail part $T_5$,  we have $T_5=P(1,2)C_{3}^{-1}$, which together with $C_{3}^{-1}=P(2,3)C_3$ implies that
$$
T_5^{-1}=P(1,2)C_3=P(1,2)P(2,3)C_3^{-1}=P(1,3)P(1,2)C_3^{-1}=P(1,3)T_5.
$$
  Finally, for the whole sequence we have $C_5^{-1}=(T_5^{-1},\ {\bf 5},\ M_5^{-1},\ {\bf 4},\ C_{3}^{-1})=P(1,2)P(4,5)C_5$ and
 $(M_5,{\bf 5},  T_5)=P(1,2)P(4,5)C_4^{-1}$.

 We call $P(1,2)P(4,5)$ the global permutation of $C_5$, and $P(1,3)$ the induced permutation of $C_5$, and denote them $P_{gl}(5)$ and $P_{in}(5)$, respectively.
Notice that  $C_4$ also has a global permutation $P_{gl}(4)=P(3,4)$ and an induced permutation $P_{in}(4)=P(1,2)$.

  For any
hyperedge $e$ of $H_{g_2(4)+g_2(3)+1}$, if it contains both color
${\bf 4}$ and color ${\bf 5}$, then it satisfies the condition of
 2-SCF-coloring. If not, then the color sequence associated with $e$
 must be a subsequence of $C_4$ or $(M_5,{\bf 5},  T_5)$ (i.e. $P(1,2)P(4,5)C_4^{-1}$), and thus
 it satisfies the condition of  2-SCF-coloring by our construction. Hence $C_5$ is a  2-SCF-coloring.

For $m=6$, following the construction of $C_5$ above, we would construct $C_6$ as follows:
\begin{eqnarray*}
&&(\underbrace{C_4,\ \ {\bf 5},\ \ T_5}_{C_5},\ \ {\bf 6},\ \ T_6),
\end{eqnarray*}
where by noticing that $T_5^{-1}=P(1,3)T_5$, we should make $T_6=P(1,3)C_4^{-1}$. Then by the fact that $C_4^{-1}=P(3,4)C_4$, we get that
$$
T_6^{-1}=P(1,3)C_4=P(1,3)P(3,4)C_4^{-1}=P(1,4)P(1,3)C_4^{-1}=P(1,4)T_6.
$$
And, we have
$$
(T_5,\ {\bf 6},\ T_6)=P(1,3)P(5,6)(C_4,\ {\bf 5},\ T_5)^{-1}=P(1,3)P(5,6)C_5^{-1},
$$
and $C_6^{-1}=(T_6^{-1},\ {\bf 6},\ T_5^{-1},\ {\bf 5},\ C_4^{-1})=P(1,3)P(5,6)C_6$. So the global permutation  $P_{gl}(6)=P(1,3)P(5,6)$ and the induced permutation  $P_{in}(6)=P(1,4)$. As to $C_5$, we can easily check that $C_6$ is a  2-SCF-coloring.

Now suppose that for any integer $m\geq 6$, we have constructed a 2-SCF-coloring $C_m$ as follows:
\begin{eqnarray*}
(\underbrace{C_{m-2},\ {\bf m-1},\ T_{m-1}}_{C_{m-1}},\ \ {\bf m}, \ \ T_m),
\end{eqnarray*}
which satisfies that
\begin{eqnarray}\label{ch-a}
T_{m-1}^{-1}=P_{in}(m-1)T_{m-1},\ T_m=P_{in}(m-1)C_{m-2}^{-1},
\end{eqnarray}
where $P_{in}(m-1)$ stands for the induced permutation of $C_{m-1}$. Denote $P_{gl}(m)=P_{in}(m-1)P(m-1,m)$. Then we have
$
C_m^{-1}=P_{gl}(m)C_m,
$
i.e. $P_{gl}(m)$ is the global permutation of $C_m$. By the second equality in (\ref{ch-a}) and the symmetry of $C_{m-2}$ (i.e. $C_{m-2}^{-1}=P_{gl}(m-2)C_{m-2}$), there exists a permutation $P_{in}(m)$ such that
$T_m^{-1}=P_{in}(m)T_m$, i.e. $P_{in}(m)$ is the induced permutation of $C_{m}$

As an algorithm we can construct $C_{m+1}$ as follows:
\begin{eqnarray}\label{ch-b}
&& (\underbrace{C_{m-1},\ {\bf m},\ T_{m}}_{C_{m}},\ \ {\bf m+1}, \ \ T_{m+1}),
\end{eqnarray}
where $T_{m+1}=P_{in}(m)C_{m-1}^{-1}$. Let $P_{gl}(m+1)=P_{in}(m)P(m,m+1)$. Then we have
\begin{eqnarray}\label{ch-c}
(T_m,\ {\bf m+1},\ T_{m+1})=P_{gl}(m+1)C_m^{-1},\ \ \  C_{m+1}^{-1}=P_{gl}(m+1)C_{m+1}.
\end{eqnarray}
By (\ref{ch-b}) and (\ref{ch-c}), we can easily check that $C_{m+1}$ is a 2-SCF-coloring and
so $g_2(m+1)\geq g_2(m)+g_2(m-1)+1$.\hfill\fbox

\begin{cor}\label{cor2.3}
For any  $n \in \mathbb{N}$, we have
\begin{eqnarray}\label{cor2.3-a}
\chi_{2scf}(H_n)=\left\lfloor\log_{\frac{1+\sqrt{5}}{2}}
\left(\sqrt{5}\left(n+\frac{1}{2}\right)\right)\right\rfloor-1,
\end{eqnarray}
where$\lfloor{ x }\rfloor$ means the largest integer that is smaller than x.
\end{cor}
{\bf Proof.} By Theorem \ref{thm2.2}, we have
$g_2(m+1)=g_2(m)+g_2(m-1)+1,\ {\forall}m\geq2$. Let
$\hat{g}_2(m)=g_2(m)+1,\ {\forall}m\geq2$.  Then $\hat{g}_2(m)$
satisfies the following recursive relation:
$$\hat{g}_2(m+1)-\hat{g}_2(m)-\hat{g}_2(m-1)=0,$$
with the initial two values $\hat{g}_2(1)=2$ and $\hat{g}_2(2)=3$.
Hence $\hat{g}_2(m)$ is just the Fibonacci number $F_{m+2}$, which can be expressed by
$$
F_{m+2}=\frac{\varphi^{m+2}-\psi^{m+2}}{\sqrt{5}},
$$
where $\varphi=\frac{1+\sqrt{5}}{2},\psi=\frac{1-\sqrt{5}}{2}$. Further, since $|\frac{\psi^{m+2}}{\sqrt{5}}|<\frac{1}{2}$, we have
$$
F_{m+2}=\left\lfloor\frac{\varphi^{m+2}}{\sqrt{5}}+\frac{1}{2}\right\rfloor.
$$
Hence we have
\begin{eqnarray}\label{cor2.3-b}
g_2(m)=\left\lfloor\frac{\varphi^{m+2}}{\sqrt{5}}+\frac{1}{2}\right\rfloor-1.
\end{eqnarray}
By (\ref{cor2.3-b}) and the formula $\chi_{2scf}(H_n)=\inf\{m:g_2(m)\geq
n\}$, we obtain (\ref{cor2.3-a}).
 \hfill\fbox

\begin{thm}\label{thm2.4}
For any integer $m \geq 3$, we have
\begin{eqnarray}\label{thm2.4-a}
g_{3}(m)=2g_{3}(m-2)+1.
\end{eqnarray}
\end{thm}
{\bf Proof.} By Corollary \ref{cor2.2}(ii), in order to prove
(\ref{thm2.4-a}), we  need only  to prove  that for any
$m\geq 3,$ $g_3(m)\geq 2g_{3}(m-2)+1$, i.e. there exists  a
3-SCF-coloring for $2g_{3}(m-2)+1$ points by using $m$ colors. As in
the proof of Theorem \ref{thm2.2}, in the following we use $(a_1,a_2,\ldots,a_n)$  to denote a
3-SCF-coloring $\phi:\{1,2,\ldots,n\}\to \{1,2,\ldots,m\}$ for some
$m$ with  $a_i=\phi(i)$ for $i=1,\ldots,n$ and denote by $C_m$ one
3-SCF-coloring for $g_{3}(m)$ points using $m$ colors.

{\bf Step 1.} Suppose $m=2p+1$ is an odd integer. For $p=1,2$, let
$C_3=({\bf 1,2,3})$ and $C_5=(1,{\bf 2},3,{\bf 4},1,{\bf 5},3)$,
where the bold colors appear only once. Denote $B_3=(1,5,3)$. Then
$C_5$ can be expressed by $(C_3, 4, B_3)$. Define $B_5=(C_3, 7, B_3)$. Then for $p=3$, we construct
$C_7$  by
\begin{eqnarray*}\label{3-C-7}
( \underbrace{C_{3},\ {\bf 4},\ B_{3}}_{C_{5}},\ {\bf 6},\
\underbrace{C_{3},\ {\bf 7},\ B_{3}}_{B_5} ).
\end{eqnarray*}
For any integer $p\geq 3$, define
\begin{eqnarray}
&&C_{2p-1}=(C_{2p-3}, 2p-2,B_{2p-3}),\ \ \ B_{2p-1}=(C_{2p-3}, 2p+1, B_{2p-3}).
\end{eqnarray}

Now assume that for $p\geq 3$ and any $l=3,\ldots,p$, we have
constructed the 3-SCF-coloring $C_{2l+1}$ by  following the above
idea.  Denote $C_{2p+1}$ by
\begin{eqnarray}\label{3-C-2p+1}
( \underbrace{C_{2p-3}, {\bf 2p-2},B_{2p-3}}_{C_{2p-1}}, \ {\bf 2p},\
\underbrace{C_{2p-3}, {\bf 2p+1}, B_{2p-3}}_{B_{2p-1}}).
\end{eqnarray}

Basing on the above coloring (\ref{3-C-2p+1}), we construct one
$(2(p+1)+1)$-coloring for $2g_3(2p+1)+1$ points as follows:
\begin{eqnarray}\label{3-C-2p+3}
 ( \underbrace{C_{2p-1}, {\bf 2p}, B_{2p-1}}_{C_{2p+1}},\ {\bf 2p+2},\
\underbrace{C_{2p-1}, {\bf 2p+3}, B_{2p-1}}_{B_{2p+1}} ).
\end{eqnarray}
Now we show that (\ref{3-C-2p+3}) is a 3-SCF-coloring. Notice that
the three colors $2p,2p+2,2p+3$ appear only once.

If a hyperedge contains the colors $2p,2p+2$ and $2p+3$, then it
satisfies the condition of  3-SCF-coloring. So we need only check
those hyperedges which do not contain all these three colors.
Hyperedges which do not contain all the colors $2p,2p+2$ and $2p+3$
have the following four types (with overlapping):
\begin{enumerate}
\item Those which do not contain color $2p+3$; \label{Non9}
\item Those which do not contain color $2p$; \label{Non6}
\item Those which do not contain color $2p+2$ and color $2p+3$; \label{Non89}
\item Those which do not contain color $2p$ and color $2p+2$. \label{Non68}
\end{enumerate}
In the following, we  only  check Type \ref{Non9} and Type
\ref{Non89}. The proofs for Type \ref{Non6} and Type \ref{Non68} are
similar to the cases of Type \ref{Non9} and Type \ref{Non89},
respectively. If a hyperedge $e$ belongs to Type \ref{Non89}, then
the color sequence associated with $e$ must be a subsequence of
$C_{2p+1}$, and thus it satisfies the condition of 3-SCF-coloring by
the inductive hypothesis.

Suppose that a hyperedge $e$ belongs to Type \ref{Non9} but not Type
\ref{Non89}. Then the color sequence associated with $e$ contains
the color $2p+2$ and is a subsequence of the following coloring
\begin{eqnarray}\label{hh}
( \underbrace{\underbrace{C_{2p-3},  2p-2,B_{2p-3}}_{C_{2p-1}}, \
{\bf 2p},\ \underbrace{C_{2p-3}, {\bf 2p+1},
B_{2p-3}}_{B_{2p-1}}}_{C_{2p+1}},{\bf 2p+2},\underbrace{C_{2p-3},
 2p-2,B_{2p-3}}_{C_{2p-1}} ).
\end{eqnarray}
Noting that in the above coloring, the three colors $2p,2p+1,2p+2$
appear only once. We know that if $e$ contains the colors $2p,2p+1$
and $2p+2$, then it satisfies the condition of  3-SCF-coloring. Thus
we need only to check that the following types of hyderedges:
\begin{enumerate}\setcounter{enumi}{4}
\item Those which  contain colors $2p+1,2p+2$ and do not contain color $2p$; \label{Non9-1}
\item Those which  contain color $2p+2$ and do not contain color $2p+1$. \label{Non6-1}
\end{enumerate}
If $e$ belongs to Type \ref{Non9-1}, then the color sequence
associated with $e$ is a subsequence of
\begin{eqnarray*}
( \underbrace{C_{2p-3}, {\bf 2p+1},B_{2p-3}}_{B_{2p-1}},{\bf
2p+2},\underbrace{C_{2p-3}, {\bf 2p-2},B_{2p-3}}_{C_{2p-1}} ),
\end{eqnarray*}
which can be obtained from $C_{2p+1}=(C_{2p-1},2p,B_{2p-1})$ by
exchanging the positions of $C_{2p-1}$ and $B_{2p-1}$, and replacing
color $2p$ with color $2p+2$. Hence in this case $e$  satisfies the
condition of  3-SCF-coloring.
 If $e$
belongs to Type \ref{Non6-1}, then the color sequence associated
with $e$ contains color $2p+2$ and is a subsequence of
$(B_{2p-3},2p+2,C_{2p-3},2p-2,B_{2p-3})$. For simplicity, we denote
the color sequence  by $D_p$. Then $D_3$ is expressed by
\begin{eqnarray*}
( \underbrace{1,5,3}_{B_3},\ {\bf 8},\ \underbrace{1,{\bf 2},3}_{C_3},\ {\bf 4},\ \underbrace{1,5,3}_{B_3} ).
\end{eqnarray*}
In this case, we can easily check that $e$  satisfies the condition
of  3-SCF-coloring. If $p\geq 4$, then $D_p$ is expressed by
\begin{eqnarray*}
( \underbrace{C_{2p-5},2p-1,B_{2p-5}}_{B_{2p-3}},{\bf
2p+2},\underbrace{C_{2p-5},{\bf 2p-4},B_{2p-5}}_{C_{2p-3}}, {\bf
2p-2},\underbrace{C_{2p-5},2p-1,B_{2p-5}}_{B_{2p-3}} ),
\end{eqnarray*}
where the three colors $2p+2,2p-4,2p-2$ appear only once as in $D_3$.
If $e$ contains color $2p-2$, then it must contain all the three colors  $2p+2,2p-4,2p-2$ and thus satisfies
the condition of  3-SCF-coloring. If $e$ does not contain color $2p-2$, then the color sequence associated
with $e$ is a subsequence of
\begin{eqnarray*}
( \underbrace{C_{2p-5},2p-1,B_{2p-5}}_{B_{2p-3}},{\bf
2p+2},\underbrace{C_{2p-5},{\bf 2p-4},B_{2p-5}}_{C_{2p-3}} ),
\end{eqnarray*}
which can be obtained from $C_{2p-1}=(C_{2p-3},2p-2,B_{2p-3})$
by exchanging the positions of $C_{2p-3}$ and $B_{2p-3}$, and replacing color $2p-2$ with color $2p+2$. Hence in this case $e$  satisfies the condition of  3-SCF-coloring.

In a word, (\ref{3-C-2p+3}) is a 3-SCF-coloring.

{\bf Step 2.} Suppose $m=2p$ is an even integer. For $p=1,2$, let
$C_2=({\bf 1,2})$ and $C_4=(1,{\bf 2},{\bf 3},{\bf 4},1)$. For
$p\geq 3$, we construct the following coloring $C_{2p}$ basing on
$C_{2p-3}$ and $B_{2p-3}$, which are defined in {\bf Step 1}, by
\begin{eqnarray}\label{even}
( \underbrace{\underbrace{C_{2p-5}, \ {2p-4},\ B_{2p-5}}_{C_{2p-3}},\
{\bf{2p-2}},\ \underbrace{C_{2p-5}, \ {\bf 2p-1},\
B_{2p-5}}_{B_{2p-3}}}_{C_{2p-1}}, \ {\bf{2p}},\
\underbrace{C_{2p-5}, \ {2p-4},\ B_{2p-5}}_{C_{2p-3}} ).
\end{eqnarray}
Denote by $|C_{2p}|$ the length of the coloring sequence  $C_{2p}$.
Then $|C_6|=3\times 3+2=11=2|C_4|+1$, and for any $p\geq 4$, we get
by (\ref{even}) and  {\bf Step 1} that
\begin{eqnarray*}
|C_{2p}|&=&3 g_3(2p-3)+2=3(2g_3(2p-5)+1)+2\\
&=&2(3g_3(2p-5)+2)+1=2|C_{2p-2}|+1.
\end{eqnarray*}
Hence it's enough to show that the coloring (\ref{even}) is  a
3-SCF-coloring. Noting that the three colors $2p-2,2p-1,2p$ appear
only once in (\ref{even}). Then following the idea in {\bf Step 1}
and by using the property of $C_{2p-1}$, we can easily obtain that
the coloring (\ref{even}) is  a 3-SCF-coloring. \hfill\fbox

\begin{cor}\label{cor2.5}
For any  $n \in \mathbb{N}$, we have
\begin{eqnarray}\label{cor2.5-a}
\chi_{3scf}(H_n)=\rm{min}\left\{
\left\lceil2(1+\rm{log}_{2}\frac{n+1}{3})\right\rceil_{e},
\left\lceil2\rm{log}_{2}(n+1)-1\right\rceil_{o}\right\},
\end{eqnarray}
where $\lceil x \rceil_{e}$ means the smallest even integer that is larger than x, and $\lceil x \rceil_{o}$ means the smallest odd integer that is larger than x.
\end{cor}
 {\bf Proof.} By Theorem \ref{thm2.4}, we have known that $g_{3}(m)=2g_{3}(m-2)+1,\
{\forall}m{\geq}3.$ Let $\hat{g}_3(m)=g_3(m)+1,\ {\forall}m\geq3$,
then $\hat{g}_3(m)$ satisfies the following recursive relation
$$\hat{g}_3(m+2)=2\hat{g}_3(m),\ \forall m\geq 1,$$
which together with  $\hat{g}_3(1)=2,\hat{g}_3(2)=3$ implies that
\begin{eqnarray}\label{cor2.5-b}
g_{3}(m)=\left\{\begin{array}{cl}
3\cdot 2^{p-1}-1 &\mbox{if}\  m=2p\ ; \\
2^{p+1}-1          &\mbox{if}\  m=2p+1.
\end{array}  \right.
\end{eqnarray}
By (\ref{cor2.5-b}) and the formula $\chi_{3scf}(H_n)=\inf\{m:g_3(m)\geq
n\}$, we obtain (\ref{cor2.5-a}).\hfill\fbox

\begin{rem}
A natural question arises:

{\bf For general $k$, the two inequalities in Corollary \ref{cor2.2}
can be strengthened to be equalities?}

\noindent We conjecture that this is true, but we can not give a
proof yet. In the following, we do some discussions.
\end{rem}

\begin{thm}\label{thm2.8}
Suppose $k,p\in\mathbb{N}$ with $k=2p+1$. Then for any  $l \in \mathbb{N}$, we have
\begin{eqnarray}\label{thm2.8-a}
g_{k}(k+l(p+1))=2g_{k}(k+(l-1)(p+1))+1.
\end{eqnarray}
\end{thm}
{\bf Proof.} By Corollary \ref{cor2.2}(ii), in order to prove
(\ref{thm2.8-a}), we  need only  to prove  that for any
$l\in \mathbb{N}$, we have $g_{k}(k+l(p+1))\geq
2g_{k}(k+(l-1)(p+1))+1$, i.e.  there exists  a $k$-SCF-coloring for
$2g_{k}(k+(l-1)(p+1))+1$ points by using $k+l(p+1)$ colors.

For simplicity, we only give the proof for $p=2$ in the following, i.e. to prove that there exists  a $5$-SCF-coloring for
$2g_5(5+(l-1)(p+1))+1$ points by using $5+l(p+1)$ colors.  For general positive integer $p$, the proof is similar.

As in
the proof of Theorem \ref{thm2.2}, in the following we use $(a_1,a_2,\ldots,a_n)$  to denote a
$5$-SCF-coloring $\phi:\{1,2,\ldots,n\}\to \{1,2,\ldots,m\}$ for
some $m$ with  $a_i=\phi(i)$ for $i=1,\ldots,n$ and denote by $C_m$
one $5$-SCF-coloring for $g_5(m)$ points using $m$ colors.

When $m=5$, let $C_5=(1,2,3,4,5)=(1,2,3,2p,2p+1).$ When $m=5+(p+1)$, we construct the
following  coloring $C_{5+(p+1)}$ basing on $C_5$ by
\begin{eqnarray}\label{thm2.8-c}
( \underbrace{1,2,3,2p,2p+1}_{C_5},
2(p+1),\underbrace{s(1),s(2),s(3),s(2p),s(2p+1)}_{B_5} ),
\end{eqnarray}
where for odd integer $i=1,3,2p+1,s(i)=i$ and for even integer $i=2j(j=1,p)$, $s(i)=2(p+1)+j$, i.e. $C_{5+(p+1)}$ can be expressed by
\begin{eqnarray}\label{thm2.8-d}
( \underbrace{1,\ {\bf 2},\ 3,\ {\bf 4},\ 5}_{C_5},\
{\bf 6},\ \underbrace{1,\ {\bf 7},\ 3,\ {\bf 8},\ 5}_{B_5} ).
\end{eqnarray}

We claim that the above coloring (\ref{thm2.8-d}) (i.e. (\ref{thm2.8-c})) is a 5-SCF-coloring. Notice that the $5$ colors in $\{2,4,6,7,8\}$ appear only once, and each color in $\{1,3,5\}$ appears twice.
Let $e$ be any hyperedge of $H_{2(2p+1)+1}$ and write $e=\{i,i+1,\ldots,i+j\}(i=1,2,\ldots,2(2p+1)+1,j=0,1,\ldots,2(2p+1)+1-i)$. If $e$ lies on the left hand side or right hand side of the point with color $2(p+1)$ i.e. color ${\bf 6}$, then obviously it satisfies the condition of $5$-SCF-coloring. Now suppose that $e$ contains the  color $2(p+1)$. If $|e|\leq 5$, then any color $i\in\{1,3,5\}$ appear at most once in $e$ since there are $5$ other colors different from $i$ between the two points with color $i$. Thus in this case, all points in $e$ have different colors.

Now we consider the case $|e|> 5$. Let $i=1$ and $j\geq 5$. If $j=5$, then the color sequence associated with $e$ is $(1,2,3,4,5,
6)$, and so all the $5+1$ colors in $\{1,2,3,4,5,6\}$ appear only once in $e$. If $j=5+q$ for $q=1,\ldots,5$, then the color sequence associated with $e$ is
$$
( 1,2,3,2p,2p+1,2(p+1),s(1),\ldots,s(q) )
$$
If $q=2l+1$ is an odd integer, then $5$ colors in the set
$
\{1,2,2p,2p+1,2(p+1)\}\cup\{s(1),\ldots,s(q)\}\backslash$ $
\{s(1),s(3),\ldots,s(q)\}$
 appear only once in $e$; if $q=2l$ is an even integer, then
$5+1$ colors in $\{1,2,\ldots,2p+1,2p+2,s(1),\ldots,s(q)\}\backslash\{s(1),s(3),\ldots,s(q-1)\}$ appear only once in $e$.
Hence $e$ satisfies the condition of 5-SCF-coloring.

For $i=2,3,\ldots,2p+1$, we can similarly show that $e$ satisfies the condition of $5$-SCF-coloring. Hence (\ref{thm2.8-d}) (i.e. (\ref{thm2.8-c})) is a 5-SCF-coloring.

Now we neglect the colors $\{1,3,5\}$ which  appear twice
in $C_{5+(p+1)}$, and only retain the colors appearing only once, to get the following  sequence of colors:
\begin{eqnarray*}
\tilde{C}_{5+(p+1)}:\ ( 2,\ 2p,\
2(p+1),\ s(2),\ s(2p)).
\end{eqnarray*}
For convenience, we write $\tilde{C}_{5+(p+1)}$ by
\begin{eqnarray*}
(\ a^{(1)}_1,a^{(1)}_2,a^{(1)}_3,a^{(1)}_{2p},a^{(1)}_{2p+1} ).
\end{eqnarray*}
Basing on $\tilde{C}_{5+(p+1)}$, we construct the following sequence of colors:
\begin{eqnarray*}
\tilde{C}_{5+2(p+1)}:\ ( \underbrace{a^{(1)}_1,a^{(1)}_2,a^{(1)}_3,a^{(1)}_{2p},a^{(1)}_{2p+1}}
_{\tilde{C}_{5+(p+1)}},
3(p+1),\underbrace{s_1(a^{(1)}_1),s_1(a^{(1)}_2),s_1(a^{(1)}_3),s_1(a^{(1)}_{2p}),
s_1(a^{(1)}_{2p+1})}
_{\tilde{B}_{5+(p+1)}} ),
\end{eqnarray*}
 where for odd integer $i=1,3,2p+1,s_1(a^{(1)}_i)=a^{(1)}_i$ and for even integer $i=2j(j=1,p)$, $s_1(a^{(1)}_i)=3(p+1)+j$, i.e. we have
 \begin{eqnarray}
\tilde{C}_{5+2(p+1)}:\ ( \underbrace{2,\ {\bf 4},\ 6,\ {\bf 7},\ 8}
_{\tilde{C}_{5+(p+1)}},\
{\bf 9},\ \underbrace{2,\ {\bf 10},\ 6,\ {\bf 11},\ 8}_{\tilde{B}_{5+(p+1)}} ).
\end{eqnarray}
 By the above sequence $\tilde{C}_{5+2(p+1)}$ and  recovering the neglected colors $\{1,3,5\}$ in $C_{5+(p+1)}$, we  construct the  following coloring $C_{5+2(p+1)}$ by
\begin{eqnarray}\label{thm2.8-e}
( C_{5+(p+1)},\ 3(p+1),\ B_{5+(p+1)} ),
\end{eqnarray}
where
\begin{eqnarray*}
C_{5+(p+1)}&=&(1,\ a^{(1)}_1,\ 3,\ a^{(1)}_2,\  5,\ a^{(1)}_3,\ 1,\ a^{(1)}_{2p},\ 3,\ a^{(1)}_{2p+1},\ 5 )\\
&=&(1,\ {\bf 2},\ 3,\ {\bf 4},\ 5,\ {\bf 6},\ 1,\ {\bf 7},\ 3,\ {\bf 8},\ 5),\\
 B_{5+(p+1)}&=&(1,\ s_1(a^{(1)}_1),\ 3,\ s_1(a^{(1)}_2),\ 5,\ s_1(a^{(1)}_3),\ 1,\ s_1(a^{(1)}_{2p}),\ 3,\
s_1(a^{(1)}_{2p+1}),\ 5)\\
&=&(1,\ {\bf 2},\ 3,\ {\bf 10},\ 5,\ {\bf 6},\ 1,\ {\bf 1}1,\ 3,\
{\bf 8},\ 5).
\end{eqnarray*}

Now we show that the coloring (\ref{thm2.8-e}) is a 5-SCF-coloring. Notice that the 5 colors in $\{a^{(1)}_2,
a^{(1)}_{2p},
3(p+1),s_1(a^{(1)}_2),
s_1(a^{(1)}_{2p})\}$ appear only once,  each color in $\{a^{(1)}_1,a^{(1)}_3,a^{(1)}_{2p+1}\}$ appears twice, and each color in $\{1,3,5\}$ appear four times.

Let $e$ be any hyperedge of $H_{2(2\times 5+1)+1}$ and write $e=\{i,i+1,\ldots,i+j\}(i=1,2,\ldots,2(2\times 5+1)+1,j=0,1,2(2\times 5+1)+1-i)$.
If $e$ lies on the left hand side or right hand side of the point with color $3(p+1)$, then it satisfies the condition of $5$-SCF-coloring by our construction and $5$-SCF property of $C_{5+(p+1)}$ and $B_{5+(p+1)}$.
Now suppose that $e$ contains the  color $3(p+1)$.
If $|e|\leq 5$, then each color in $\{1,3,5\}\cup \{a^{(1)}_1,a^{(1)}_3,a^{(1)}_{2p+1}\}$ appear at most once in $e$ since for any color $i\in \{1,3,5\}$, there are at least $5$ other colors different from $i$ between any two points with color $i$; and for any color $i\in \{a^{(1)}_1,a^{(1)}_3,a^{(1)}_{2p+1}\}$, there are $2\times 5+1$ other colors different from $i$ between the two points with color $i$. Thus in this case, all the colors in $e$ are different.

Now we consider the case $|e|>5$. Let $i=1$ and $j\geq 5$. If
$5\leq j\leq 2\times 5+1$, then the color sequence associated with $e$ is a
subsequence of $(C_{5+(p+1)},3(p+1))$, and thus  it satisfies the
condition of $5$-SCF-coloring. If $j=2\times 5+2,\ldots,2(2\times 5+1)$, then all
the colors in $\{1,3,5\}$ appear at least twice in $e$.
Thus $e$ satisfies the condition of $5$-SCF-coloring if and only if
there exist at least $5$ colors in $\tilde{C}_{5+2(p+1)}$ such that
they appear only once in $e$, which is true by the fact that the
color sequence associated with $e$ contains
$(\tilde{C}_{5+(p+1)},3(p+1))$ as a subsequence, and the proof that
(\ref{thm2.8-c}) is a $5$-SCF-coloring.

For $i=2,3,\ldots,2\times 5+1$, we can similarly show that $e$ satisfies the condition of $5$-SCF-coloring. Hence (\ref{thm2.8-e}) is a $5$-SCF-coloring.

The above constructions for two colorings $C_{5+(p+1)}$ and
$C_{5+2(p+1)}$ are recursive. Hence for any  $l\in \mathbb{N}$, there
exists  a $5$-SCF-coloring for $2g_{5}(5+(l-1)(p+1))+1$ points by
using $5+l(p+1)$ colors. The proof is complete. \hfill\fbox

\begin{cor}\label{cor2.8}
Suppose $k,p\in\mathbb{N}$ with $k=2p+1$. Then for any  $l \in \mathbb{N}$, we have
\begin{eqnarray}\label{cor2.8-a}
g_{k}(k+l(p+1))=2^l(k+1)-1.
\end{eqnarray}
\end{cor}
{\bf Proof.} For any $l \in \mathbb{N}$, let $\hat{g}_l:=g_{k}(k+l(p+1))+1$. Then by Theorem \ref{thm2.8}, we have
$$
\hat{g}_l=2\hat{g}_{l-1},\ \forall l\in \mathbb{N},
$$
which together with $\hat{g}_0=k+1$ implies that
$\hat{g}_l=2^l(k+1)$. Hence (\ref{cor2.8-a}) holds.\hfill\fbox

\bigskip

We remark  that Horev et al.
\cite{22} show that $\chi_{kscf}(H_n)\leq k\log n$ (as a special  case of a more general framework), and Gargano and Rescigno \cite{GR12} show that
$$
\left\lceil\frac{k}{2}\right\rceil\left\lceil\log_2\frac{n}{k}\right\rceil
\leq \chi_{kscf}(H_n)\leq k(\left\lfloor\log_2\left\lceil\frac{n}{k}\right\rceil\right\rfloor+1).
$$
By Corollary \ref{cor2.8}, we can get the following bounds.

\begin{cor}\label{cor2.9}
(i) If $k=2p+1$ is an odd natural number. Then for any integer  $n\geq k+1,$ we have
$$
k+\frac{k+1}{2}\left(\log_2\frac{n+1}{k+1}-1\right)< \chi_{kscf}(H_n)< k+\frac{k+1}{2}\left(1+\log_2\frac{n+1}{k+1}\right).
$$
(ii) If $k=2p$ is an even natural number. Then for any integer $n\geq k+2,$ we have
$$
k-1+\frac{k}{2}\left(\log_2\frac{n+1}{k}-1\right)< \chi_{kscf}(H_n)< 1+k+\frac{k+2}{2}\left(1+\log_2\frac{n+1}{k+2}\right).
$$
\end{cor}
{\bf Proof.} (i) Let $2^{l-1}(k+1)-1<n\leq 2^l(k+1)-1$ for some integer $l\geq 1,$ then we have
\begin{eqnarray*}
\log_2\frac{n+1}{k+1}\leq l<1+\log_2\frac{n+1}{k+1}.
\end{eqnarray*}
Hence by Corollary \ref{cor2.8}, we have
\begin{eqnarray*}
\chi_{kscf}(H_n)\leq k+l(p+1)&< & k+\left(\frac{k-1}{2}+1\right)\left(1+\log_2\frac{n+1}{k+1}\right)\\
&=&k+\frac{k+1}{2}\left(1+\log_2\frac{n+1}{k+1}\right),
\end{eqnarray*}
and
\begin{eqnarray*}
\chi_{kscf}(H_n)> k+(l-1)(p+1)&\geq & k+\left(\frac{k-1}{2}+1\right)\left(\log_2\frac{n+1}{k+1}-1\right)\\
&=&k+\frac{k+1}{2}\left(\log_2\frac{n+1}{k+1}-1\right).
\end{eqnarray*}
(ii) By the monotone property of $\chi_{kscf}(H_n)$ with respect to $k$ and (i), we have
\begin{eqnarray*}
\chi_{kscf}(H_n)\leq \chi_{(k+1)scf}(H_n)< 1+k+\frac{k+2}{2}\left(1+\log_2\frac{n+1}{k+2}\right),
\end{eqnarray*}
and
\begin{eqnarray*}
\chi_{kscf}(H_n)\geq \chi_{(k-1)scf}(H_n)> k-1+\frac{k}{2}\left(\log_2\frac{n+1}{k}-1\right).
\end{eqnarray*}
\hfill\fbox

\section{$k$-CF-coloring of $H_n$}\setcounter{equation}{0}\label{kCF}
In this section, we consider $k$-CF-coloring of $H_n$ for any
 $k\in \mathbb{N}$, and obtain the following result.

\begin{thm}\label{thm3.1}
 For any  $k,n\in \mathbb{N}$, we have $\chi_{kcf}(H_n)=
 \lfloor log_{(k+1)}n \rfloor + 1$.
\end{thm}
{\bf Proof.} {\bf Step 1:} We prove for any  $m\in \mathbb{N}$ when
$n\geq (k+1)^{m}$, $\chi_{kcf}(H_n)\geq m+1$. If $m=1$, it's true.
Suppose that the claim holds for some $m\in \mathbb{N}$. We will
show that the claim holds for $m+1$ by the inductive method, i.e.
for
  $n\geq (k+1)^{m+1}$,  we will prove that
  $\chi_{kcf}(H_n)\geq m+2$.

Let $n\geq (k+1)^{m+1}$. Express $V_n=\{1,2,\ldots,n\}$ by
$V_n=V_{n,1}\cup V_{n,2}\cup \cdots\cup V_{n,k+1}\cup V_{n,k+2},$
where
\begin{eqnarray*}
&&V_{n,1}=\{1,2,\ldots,(k+1)^m\},\\
&&V_{n,2}=\{(k+1)^m+1,(k+1)^m+2,\ldots,2(k+1)^m\},\\
&&\quad\quad\quad\cdots\\
&&V_{n,k+1}=\{k(k+1)^m+1,k(k+1)^m+2,\ldots,(k+1)^{m+1}\},
\end{eqnarray*}
and if $n=(k+1)^{m+1}$, then $V_{n,k+2}=\emptyset$; if
$n>(k+1)^{m+1}$, then $V_{n,k+2}=\{(k+1)^{m+1}+1,\ldots,n\}$.
Now
suppose that we have a coloring $\phi: V_n\to \{1,2,\ldots,l\}$ for
some integer $l$ such that it is a $k$-CF-coloring for $H_n$. For
any $i=1,\ldots,k+1$, denote by  $|\phi(V_{n,i})|$ the number of
colors in $\phi(V_{n,i})$. Then by the inductive hypothesis,
$|\phi(V_{n,i})|\geq m+1$. Hence
$|\phi(\cup_{i=1}^{k+1}V_{n,i})|\geq m+2$, because otherwise all the
$m+1$ colors in $\phi(\cup_{i=1}^{k+1}V_{n,i})$ appear at least
$k+1$ times and thus the coloring $\phi$ is not a $k$-CF-coloring.
Hence $\chi_{kcf}(H_n)\geq m+2$.

{\bf Step 2:} We prove that for   $m\in \mathbb{N},n=(k+1)^m-1$,
$\chi_{kcf}(H_n)=m$. When $m=1$, $n=k$ and so $\chi_{kcf}(H_n)=1=m$.
Suppose that the claim holds for some $m\in\{1,2,\ldots\}$. For $n=
(k+1)^{m+1}-1$, we express $V_n=\{1,2,\ldots,n\}$ by $
V_n=\overline{V}_{n,1}\cup \{(k+1)^m\}\cup \overline{V}_{n,2}\cup
\{2(k+1)^m\}\cup\cdots\cup \overline{V}_{n,k}\cup\{k(k+1)^m\}\cup
\overline{V}_{n,k+1}, $ where
\begin{eqnarray*}
&&\overline{V}_{n,1}=\{1,2,\ldots,(k+1)^m-1\},\\
&&\overline{V}_{n,2}=\{(k+1)^m+1,(k+1)^m+2,\ldots,2(k+1)^m-1\},\\
&&\quad\quad\quad\cdots\\
&&\overline{V}_{n,k+1}=\{k(k+1)^m+1,k(k+1)^m+2,\ldots,(k+1)^{m+1}-1\}.
\end{eqnarray*}
For any $i=1,\ldots,k+1$, by the inductive hypothesis, we  know that
the sub-hypergraph induced by $\overline{V}_{n,i}$ has a
$k$-CF-coloring by using $m$ colors e.g. colors $1,2,\ldots,m$. Use
color $m+1$ to color the vertices
$(k+1)^m,2(k+1)^m,\ldots,k(k+1)^m$. We can easily check that this
coloring is a $k$-CF-coloring of $H_n$. So $\chi_{kcf}(H_n)\leq
m+1$. By {\bf Step 1} and the fact that $(k+1)^{m+1}-1\geq
(k+1)^{m}$, we get that $\chi_{kcf}(H_n)=m+1$ for $n=(k+1)^{m+1}-1$.
Hence for any  $m\in \mathbb{N}$, the claim holds.

By {\bf Step 1} and {\bf Step 2}, we obtain that for any
$k, n\in \mathbb{N}$, we have $\chi_{kcf}(H_n)=
 \lfloor log_{(k+1)}n \rfloor + 1$.\hfill\fbox

\newpage

{ \noindent {\bf\large Acknowledgments}

The authors acknowledge the helpful suggestions and comments of three
anonymous referees, which helped improve the first three versions  of this manuscript, respectively. Research supported
by NNSFC, Jiangsu Province basic research program (Natural Science
Foundation) (Grant No. BK2012720).

\bibliographystyle{amsplain}

\end{document}